\documentclass[preprint,12pt]{elsarticle}
\makeatletter
\def\ps@pprintTitle{%
 \let\@oddhead\@empty
 \let\@evenhead\@empty
 \def\@oddfoot{\centerline{\thepage}}%
 \let\@evenfoot\@oddfoot}
\makeatother
\usepackage{graphicx}   % if compiled with pdflatex
\usepackage{amsmath,amsfonts,amssymb,amsthm}

\usepackage[colorlinks=true,linkcolor=blue]{hyperref}

\usepackage{txfonts}%
\usepackage{bm}

\newtheorem{remark}{Remark}
\newtheorem{proposition}{Proposition}

\newcommand{\rank}{\operatorname{rank}}
\newcommand{\tr}{\operatorname{tr}}
\newcommand{\hFoo}{{\kern.15em{}_{1}\kern-.05em {F}_{1}}}

\newcommand{\mb}{\mathbf}
\newcommand{\ua}{{\mb{a}}}
\newcommand{\ub}{{\mb{b}}}
\newcommand{\um}{{\mb{m}}}
\newcommand{\uw}{{\mb{w}}}
\newcommand{\ux}{{\mb{x}}}
\newcommand{\uA}{{\mb{A}}}
\newcommand{\uH}{{\mb{H}}}
\newcommand{\uM}{{\mb{M}}}
\newcommand{\uP}{{\mb{P}}}
\newcommand{\uQ}{{\mb{Q}}}
\newcommand{\uT}{{\mb{T}}}
\newcommand{\uW}{{\mb{W}}}
\newcommand{\uX}{{\mb{X}}}
\newcommand{\usigma}{{\pmb{\sigma}}}
\newcommand{\uSigma}{{\pmb{\Sigma}}}
\newcommand{\uOmega}{{\pmb{\Omega}}}
\newcommand{\mzero}{\mb{0}}
\newcommand{\mbI}{\mb{I}}

\begin{document}

\begin{frontmatter}
%\title{Hypergeometric series and weights of chi-square distributions in a Schur complement in a noncentral Wishart matrix}
%\title{Mixture Weights of Chi-Square Distributions \\ for the Scalar Schur Complement in a \\ Noncentral Wishart Matrix}
\title{Chi-Square Mixture Representations for the Distribution of the Scalar Schur Complement in a Noncentral Wishart Matrix}

\author{Constantin Siriteanu\fnref{label_Costi}}
\ead{constantin.siriteanu@ist.osaka-u.ac.jp}
\address[label_Costi]{Graduate School of Information Science and Technology, Osaka University, Osaka, Japan}
\author{Satoshi Kuriki\fnref{label_Kuriki}}
\ead{kuriki@ism.ac.jp}
\address[label_Kuriki]{The Institute of Statistical Mathematics, Tokyo, Japan}
\author{Donald Richards\fnref{label_Donald}}
\ead{richards@stat.psu.edu}
\address[label_Donald]{Department of Statistics, Pennsylvania State University, University Park, Pennsylvania, USA}
\author{Akimichi Takemura\fnref{label_Takemura}}
\ead{takemura@stat.t.u-tokyo.ac.jp}
\address[label_Takemura]{Graduate School of Information Science and Technology, University of Tokyo, Tokyo, Japan}

\begin{abstract}
%We show that the distribution of the scalar Schur complement in a noncentral Wishart matrix is a mixture of noncentral chi-square distributions and
%hence can also be written as a mixture of central chi-square distributions 
%with different degrees of freedom.
%In deriving these results, it is crucial to express the weights of the mixture representation in terms of known discrete probability distributions over the set of nonnegative integers;
%for the case of a rank-1 noncentrality matrix, we prove that 
%the weights are the probabilities arising from a noncentral beta mixture of Poisson distributions.
We show that the distribution of the scalar
Schur complement in a noncentral Wishart matrix is 
a mixture of central chi-square distributions with different degrees of freedom.
For the case of a rank-1 noncentrality matrix, the weights of the mixture representation
arise from a noncentral beta mixture of Poisson distributions.
\end{abstract}

\begin{keyword}
Confluent Appell function; noncentral beta distribution.
\end{keyword}

\end{frontmatter}

\section{Introduction}
\label{sec:intro}

\subsection{Motivating the study of the Schur complement}

This paper considers the distribution of the Schur complement in a real-valued noncentral Wishart matrix. 
The matrix Schur complement arises naturally in many areas of multivariate statistical analysis as the covariance matrix of the conditional Gaussian distribution, in decompositions of the Wishart matrix, and in hypothesis testing problems \cite{anderson-3rd}, \cite{cottle_laa_74}, \cite[Ch.~6]{zhang_book_05}.
In particular, the Schur complement is a crucial quantity in multivariate statistical inference with monotone incomplete data \cite{chang-richards-2009}, \cite{chang-richards-2010}.

The Schur complement also appears in multiple-input/multiple-output (MIMO) wireless communications systems, as follows. 
For a communications link over a randomly changing radio channel, the signal-to-noise ratio is a crucial metric that determines all other performance measures, e.g., error probability, outage probability, capacity. 
The signal-to-noise ratio is proportional with the scalar Schur complement in the sample correlation matrix formed with the random MIMO channel matrix when the received data is estimated by the least-squares-like method known as {\it zero-forcing detection} \cite{siriteanu_twc_13}, \cite{siriteanu_twc_SC_14}, \cite{siriteanu_twc_15}.
Knowing the distribution of the signal-to-noise ratio allows us to characterize the performance measures statistically, e.g., in average.  

In the case of MIMO channel matrices that are complex-valued Gaussian distributed, we have derived recently analyses of performance measures by means of the Schur complement. For instance, in \cite{siriteanu_twc_13}, \cite{siriteanu_twc_15} we characterized with infinite series the distribution of the scalar Schur complement in a complex-valued Wishart matrix with rank-1 noncentrality matrix.  Furthermore, in \cite{siriteanu_twc_SC_14}, we characterized the more general matrix Schur complement, under certain conditions. 
However, the weights in those infinite series \cite[Eq.~(39)]{siriteanu_twc_13}, \cite[Eq.~(59)]{siriteanu_twc_15} are alternately positive and negative, which can lead to numerical instability in evaluating ensuing MIMO performance measures \cite[Sec. V.F]{siriteanu_twc_13}, and that phenomenon partly motivated us to research mixture representations for the distribution of the scalar Schur complement, which we now do in this paper.

\subsection{Contributions and Organization of the Paper}

In this paper we consider the case of real-valued Gaussian distributed matrix with rank-1 mean, 
which yields a real-valued Wishart matrix with rank-1 noncentrality matrix.
We first show that the distribution of the scalar Schur complement in this Wishart matrix is 
a mixture of noncentral $\chi^2$ distributions, where the noncentrality parameter is random.
Since the noncentral $\chi^2$ distribution itself is a Poisson mixture of central $\chi^2$
distributions with different degrees of freedom, 
the distribution of the scalar Schur complement is also
a mixture of central $\chi^2$ distributions with different degrees of freedom.

We also express the mixture weights in terms of known discrete
probability distributions over the set of nonnegative integers.
Whereas, in general, these weights are difficult to evaluate,
for the case of a rank-1 noncentrality matrix, we prove that
the weights are the probabilities arising from a noncentral beta mixture of Poisson distributions.

Thus, we prove in Section \ref{sec:sc} that the distribution of the scalar Schur complement is a mixture of noncentral $\chi^2$ distributions, hence also a mixture of central $\chi^2$ distributions with different degrees of freedom.
In Section \ref{sec:pdf} we relate the moment-generating function of the scalar Schur complement to the probability-generating function of the mixture weights of central $\chi^2$ distributions.
In Section \ref{sec:rank1} we derive the distribution of the scalar Schur complement in a noncentral Wishart matrix with a rank-1 noncentrality matrix.

\subsection{Notation}

\begin{itemize}
\item Scalars, vectors, and matrices are represented with lowercase italics, lowercase boldface, and uppercase boldface, respectively, e.g., $u$, $\ux$, and $\uX$; %the notation $\uX \in \mathbb{R}^{n \times p}$ indicates $n$ rows and $p$ columns for $\uX$; 
$\mathbb{R}^{n}$ and $\mathbb{R}^{n\times p}$ denote the set of $n$-dimensional column vectors and the set of $n\times p$ matrices, respectively;
The zero vectors and matrices of appropriate dimensions are denoted by $\mzero$; the superscript $\cdot'$ stands for transpose;
$\mbI_n$ denotes the $n \times n$ identity matrix.
%\item $[\cdot]_{i}$ is the $i$th element of a vector; $[\cdot]_{i,j}$, $[\cdot]_{i,\bullet}$, and $[\cdot]_{\bullet,j}$ indicate the $i,j$th element, $i$th row, and $j$th column of a matrix; $\| \uH \|^2 = \sum_{i=1}^{\NR} \sum_{j=1}^{\NT} |[\uH]_{i,j}|^2$ is the squared Frobenius norm.
%\item $i = 1:N$ stands for the enumeration $i = 1, \, 2, \, \ldots, \, N$.
\item The notation $\ux \sim {\cal{N}}_n \left(\um, \uSigma \right)$ indicates an $n$-dimensional, real-valued, Gaussian random vector with mean $\um$ and covariance $\uSigma$;
$\uW \sim {\cal{W}}_p \left(n, \uSigma, \Omega \right)$ indicates that $ \uW $ has a $p\times p$ noncentral Wishart distribution with $n$ degrees of freedom, covariance matrix $\uSigma$, and noncentrality matrix $\uOmega$ \cite[Eq.~(67)]{james_ams_64};
$\ux_1 \,|\, \uX_2$ stands for $\ux_1$ conditioned on $\uX_2$; $\mathbb{E}[\cdot]$ denotes expectation;
$\chi_\nu^2(\delta)$ denotes the noncentral \textit{chi-square} distribution with $\nu$ degrees of freedom and noncentrality parameter $\delta$;
$\chi_\nu^2$ denotes the central chi-square distribution with $\nu$ degrees of freedom;
$Beta(k,l)$ represents the central beta distribution with shape parameters $k$ and $l$;
$Beta(k,l,\delta)$ represents the noncentral beta distribution with shape parameters $k$ and $l$, and noncentrality $\delta$;
finally, $Pois(\lambda)$ represents the Poisson distribution with mean $\lambda$.
\item
For ${\rm Re}(k), {\rm Re}(l) > 0$, ${\Gamma}(k) = \int_{0}^{\infty} t^{k-1} e^{-t} {\rm d}t$ denotes the gamma function and 
$B(k,l) = \int_{0}^{1} t^{k-1} (1-t)^{l-1} {\rm d}t$ denotes the beta function;
$(x)_n$ is the rising factorial, where $(x)_0 = 1$ and $(x)_n = x (x + 1) \cdots (x + n - 1)$, $n \geq 1$;
$\hFoo (\cdot; \cdot; \cdot)$ is the {confluent hypergeometric function}~\cite[Eq.~(13.2.2), p.~322]{NIST_book_10}:
\[
 \hFoo(b;c;w) = \sum_{l=0}^\infty \frac{(b)_l}{(c)_l}\frac{w^l}{l!};
\]
and $\Phi_2(\cdot,\cdot;\cdot;\cdot,\cdot)$ is the confluent Appell function
\cite[Eq.~(1.3)]{brychkov-saad-2012}:
\[
 \Phi_2(b,b';c;w,z) = \sum_{l=0}^\infty\sum_{m=0}^\infty \frac{(b)_l(b')_m}{(c)_{l+m}}\frac{w^l z^m}{l! m!}.
\]
\item Acronyms and abbreviations: 
 %GHFD for generalized hypergeometric factorial moment distribution,
 p.d.f.~for the probability density function, p.g.f.~for the probability generating function, m.g.f.~for the moment-generating function.
\end{itemize}
%
%\subsection{Organization}
%
%We prove in Section \ref{sec:sc} that the distribution of the scalar Schur complement is a mixture of noncentral $\chi^2$ distributions, hence also a mixture of central $\chi^2$ distributions with different degrees of freedom.
%In Section \ref{sec:pdf} we relate the moment-generating function of the scalar Schur complement to the probability-generating function of the mixture weights of central $\chi^2$ distributions.
%In Section \ref{sec:rank1} we derive the distribution of the scalar Schur complement in a noncentral Wishart matrix with a rank-1 noncentrality matrix.

\section{Distribution of the scalar Schur complement in a noncentral Wishart matrix}
\label{sec:sc}

\subsection{The scalar Schur complement in a noncentral Wishart matrix}

Let $\uX = (x_{ij}) \in \mathbb{R}^{n \times p}$ be a real-valued Gaussian random matrix with mean $\uM = (\mu_{ij}) \in \mathbb{R}^{n \times p}$, and whose rows are mutually independent $p$-dimensional vectors with covariance matrix $\uSigma = (\sigma_{ij}) \in \mathbb{R}^{p \times p}$. 
%Note that then we can write
%\begin{eqnarray}
%\label{equation_Sigma_X_M}
%\uSigma = \frac{1}{n} \mathbb{E}[ ( \uX - \uM )' (\uX -\uM) ].
%\end{eqnarray}

Let us partition $\uX$ and $\uM$ into the first column and the rest of columns, i.e.,
\[
\uX=\begin{pmatrix} \ux_1 & \uX_2\end{pmatrix}, \qquad \uM=\begin{pmatrix}\um_1 & \uM_2 \end{pmatrix},
\]
with $\ux_1, \um_1 \in \mathbb{R}^{n}$, and $\uX_2, \uM_2 \in \mathbb{R}^{n \times (p-1)}$. 
Accordingly, we partition $\uSigma$ as
\begin{eqnarray*}
%\label{equation_covariance}
\uSigma=\begin{pmatrix} \sigma_{11} & \usigma_{21}' \\ \usigma_{21} & \uSigma_{22} \end{pmatrix}.
\end{eqnarray*}
Then, the conditional distribution of $\ux_1$ given $\uX_2$ is~\cite[p.~35]{muirhead_book_05}
\begin{equation}
\label{eq:x1givenx2}
\ux_1 \,|\, \uX_2 \sim {\cal{N}}_n(\um_1+ (\uX_2 -\uM_2 ) \uSigma_{22}^{-1} \usigma_{21}, \, \sigma_{11\cdot 2} \mbI_n),
\end{equation}
where 
\begin{equation}
\label{equation_covariance}
\sigma_{11\cdot 2} = \sigma_{11} - \usigma_{21}' \uSigma_{22}^{-1} \usigma_{21}
\end{equation}
is the variance of each element of $\ux_1 \,|\, \uX_2$, and, by definition, also the Schur complement of $\uSigma_{22}$ in $\uSigma$.

Now, let us consider the matrix of interest, 
% , the joint probability density function (p.d.f.) of the elements of $\uX$ is given by~\cite[Eq.~(1)]{james_ams_64}.
% Further, the sample correlation matrix, i.e.,
\[
\uW= \uX' \uX \in \mathbb{R}^{p \times p}.
\]
%which, based on~(\ref{equation_Sigma_X_M}), satisfies
%\begin{eqnarray}
%\uSigma = \frac{1}{n} \mathbb{E}[ \uW ] - \frac{1}{n} \mathbb{E}[ \uM' \uM ].
%\end{eqnarray}
%Thus, we may regard $\uW$ as the sample correlation matrix of the rows of $\uX$.
It is well known that for $\uM \neq \mzero$ we have $\uW \sim {\cal{W}}_p \left(n, \uSigma, \uOmega \right)$ with $\uOmega=\uSigma^{-1}\uM' \uM$.
%Because, based on~(\ref{equation_Sigma_X_M}), we have that
%\begin{eqnarray}
%\uSigma = \frac{1}{n} \mathbb{E}[ \uW ] - \frac{1}{n} \mathbb{E}[ \uM' \uM ],
%\end{eqnarray}
%we can view $\uW$ as the sample correlation matrix of the rows of $\uX$.
%
If we partition $\uW$ analogously to~(\ref{equation_covariance}) as
\[
%\label{equation_partitioned_RTK}
\uW = \begin{pmatrix}
    w_{11} & \uw_{21}' \\
    \uw_{21} & \uW_{22}
  \end{pmatrix} =
  \begin{pmatrix}
    \ux_1'\ux_1 & \ux_1' \uX_2 \\
    \uX_2' \ux_1 & \uX_2' \uX_2
  \end{pmatrix}
\]
then, by definition, the Schur complement of $\uW_{22}$ in $\uW$, which is the scalar Schur complement of interest here, is given by
\begin{align}
%\label{equation_SC_definition}
w_{11\cdot 2}
&= w_{11} - \uw_{21}' \uW_{22} \uw_{21} \nonumber \\
%\label{equation_SC_sample_correlation}
&= \ux_1' \ux_1 - \ux_1'  \uX_2 (\uX_2'\uX_2)^{-1} \uX_2' \ux_1 \nonumber \\
&= \ux_1' (\mbI_n - \uX_2 (\uX_2'\uX_2)^{-1} \uX_2') \ux_1 \nonumber \\
%\label{equation_SC_definition_Hermitian_form}
&= \ux_1' \uQ_2 \ux_1, 
\label{equation_SC_definition}
\end{align}
where
\[
%\label{equation_Q2}
\uQ_2 = \mbI_n - \uX_2 (\uX_2'\uX_2)^{-1} \uX_2'.
\]

\subsection{Distribution of the conditioned scalar Schur complement}

%Consider the conditional distribution of $\ux_1$ given $\uX_2$. Since the rows of $\uX$ are independent,
%given $\uX_2$, each element $x_{i1}$ of $\ux_1$ is independent and has the
%Gaussian distribution~\cite[p.~35]{muirhead_book_05}
%\begin{eqnarray}
%{\cal{N}} \left( \mu_{i1} + \left( (x_{i2},\dots,x_{ip}) - (\mu_{i2},\dots,\mu_{ip}) \right) \uSigma_{22}^{-1} \usigma_{21}, \,
%\sigma_{11\cdot 2} \right). 
%\end{eqnarray}
%Hence, the conditional distribution of $\ux_1$ given $\uX_2$ is 
%\begin{equation}
%\label{eq:x1givenx2}
%\ux_1 \,|\, \uX_2 \sim {\cal{N}}_n( \tilde \um_1 + \uX_2 \uSigma_{22}^{-1} \usigma_{21}, \, \sigma_{11\cdot 2} \mbI_n), \qquad 
%\tilde \um_1 = \um_1 -\uM_2 \uSigma_{22}^{-1} \usigma_{21}. 
%\end{equation}
%Consider the conditional distribution of $\ux_1$ given $\uX_2$, i.e.,~\cite[p.~35]{muirhead_book_05}
%\begin{equation}
%\label{eq:x1givenx2}
%\ux_1 \,|\, \uX_2 \sim {\cal{N}}_n(\tilde \um_1 + \uX_2 \uSigma_{22}^{-1} \usigma_{21}, \, \sigma_{11\cdot 2} \mbI_n), \qquad
%\tilde \um_1 = \um_1 -\uM_2 \uSigma_{22}^{-1} \usigma_{21}. 
%\end{equation}
Based on $\ux_1 \,|\, \uX_2$ characterized in~(\ref{eq:x1givenx2}) we have
\begin{equation}
\label{eq:x1-sub-X2}
\tilde \ux_1 = \ux_1 - \uX_2 \uSigma_{22}^{-1} \usigma_{21} \sim {\cal{N}}_n(
\tilde \um_1, \sigma_{11\cdot 2} \mbI_n),
\end{equation}
where
\[
 \tilde \um_1 = \um_1 -\uM_2 \uSigma_{22}^{-1} \usigma_{21}.
\]
Note that $\tilde \ux_1$ is independent of $\uX_2$, and that the distribution of $\ux_1 \,|\, \uX_2$ is the distribution of $
 \tilde \ux_1 + \uX_2 \uSigma_{22}^{-1} \usigma_{21} $.
By substituting $\ux_1$ from~(\ref{eq:x1-sub-X2}) into~(\ref{equation_SC_definition})
 and using $\uX_2' \left (\mbI_n - \uX_2 (\uX_2'\uX_2)^{-1} \uX_2' \right) \uX_2 = \mzero$, we can write the Schur complement of interest as
\[
w_{11\cdot 2} =\tilde \ux_1' \left(\mbI_n - \uX_2 (\uX_2'\uX_2)^{-1} \uX_2' \right) \tilde \ux_1.
\]
%It follows then that
%\begin{eqnarray}
%\frac{w_{11\cdot 2}}{\sigma_{11\cdot 2}} \bigg| \uX_2 \sim \chi^2_{\nu=n-p+1} \left( {\delta} = \frac{\tilde \um_1' (\mbI_n - \uX_2 (\uX_2'\uX_2)^{-1} \uX_2) \tilde \um_1}{\sigma_{11\cdot 2}}\right).
%\end{eqnarray}

Now, let us define:
\begin{align}
%\label{equation_definition_rho1}
{\rho} &= \frac{w_{11\cdot 2}}{\sigma_{11\cdot 2}} = \frac{\tilde \ux_1' \uQ_2 \tilde \ux_1}{\sigma_{11\cdot 2}}, \nonumber \\
\label{equation_definition_u_m1_Q2}
u &= \frac{\tilde \um_1' \uQ_2 \tilde \um_1}{\Vert \tilde \um_1\Vert^2}, \\
\label{equation_definition_lambda}
\lambda &= \frac{\Vert \tilde \um_1\Vert^2}{\sigma_{11\cdot 2}}, \\
\label{equation_definition_delta1}
{\delta} &= \lambda u =\frac{\tilde \um_1' \uQ_2 \tilde \um_1}{\sigma_{11\cdot 2}}.
\end{align}
%We call ${\rho}$ the {\it scaled scalar SC}.

Because the $n \times n$ matrix $\uQ_2$ is idempotent, i.e., the magnitude of its eigenvalues is either $0$ or $1$, we deduce that $u \in [0,1]$. 
Furthermore, the fact that $\uQ_2$ has $p - 1$ eigenvalues equal to $0$ and $\nu=n - p + 1$ eigenvalues equal to $1$, yields, in general, %
%\footnote{Herein, we assume that $\tilde \um_1\neq 0$ because the interesting case with $\tilde \um_1= \mzero$, whereby ${\rho} \sim \chi^2_{\nu}$ irrespective of the distribution of $\uX_2$, has been studied in~\cite{siriteanu_twc_SC_14}.}
the following noncentral $\chi^2$ distribution for the scaled Schur complement ${\rho}$ conditioned on $u$~\cite{mathai_book_95}:
\begin{equation}
\label{equation_SC_cond_distrib_u}
{\rho} \,|\, u \sim \chi^2_{\nu } ({\delta}), \quad \nu=n-p+1.
\end{equation}
Thus, the p.d.f. of ${\rho} \,|\, u$ is given by \cite[p.~27]{mathai_book_95} as
\begin{equation}
\label{eq:non-central-chi-square-density}
g(w;\nu, {\delta}) = \sum_{k=0}^\infty \frac{({\delta}/2)^k}{k!}e^{- {\delta}/2} g_{\nu+2k}(w), \ \ w > 0,
\end{equation}
where
\[
%\label{eq:non-central-chi-square-density_coeff}
g_m(w) = \frac{1}{2^{m/2}{\Gamma}(m/2)} w^{m/2-1} e^{-w/2}
\]
is the p.d.f.~of a $\chi^2_{m}$ random variable. 

\section{Relationship between the mixture m.g.f.~and the mixture-weights p.g.f.}
% of central chi-square distributions with different degrees of freedom}
\label{sec:pdf}

\subsection{M.g.f.--p.g.f.~relationship}

Note first that the distribution of ${\rho} \,|\, u$, characterized by the p.d.f.~$g(w;\nu, {\delta})$ in~(\ref{eq:non-central-chi-square-density}), is a mixture of central $\chi^2$ distributions. 
Note also that the mixture weight, i.e., 
\begin{equation}
\label{equation_new_alpha_k}
\alpha_k = \frac{({\delta}/2)^k}{k!}e^{- {\delta}/2}
\end{equation}
represents the probability that the Poisson random variable $A \sim Pois({\delta}/2)$ takes the value $k$.  
On the one hand, the p.g.f.~of $A$ is given by
\begin{equation}
\label{equation_pgf_GA_s}
 G_A(s) = \mathbb{E}\left[s^A\right] = \sum_{k=0}^\infty s^k \alpha_k = e^{({\delta}/2)(s-1)}.
\end{equation}
On the other hand, the m.g.f.~of ${\rho} \,|\, u$ is given by the expression
\begin{equation}
\label{eq:non-central-chi-mgf}
M_{{\rho} \,|\, u }(\theta) = \mathbb{E}\left[e^{\theta {\rho}} \,|\, u\right] 
= (1-2\theta)^{-\nu/2} \exp\left(\frac{{\delta}\theta}{1-2\theta}\right),
\end{equation}
which is simpler than the corresponding p.d.f.~expression in~(\ref{eq:non-central-chi-square-density}).
%It is often easier to manipulate the m.g.f.~instead of the p.d.f., as also noticeable in this case by comparing~(\ref{eq:non-central-chi-square-density}) and~(\ref{eq:non-central-chi-mgf}).
%Further, by definition, the p.g.f.~of $N$ is given by
%\begin{eqnarray}
%G_N(s)=\mathbb{E}\left[s^N\right] = \sum_{k=0}^\infty s^k \beta_k.
%\end{eqnarray}
%\subsection{Distribution of the unconditioned SC}

\begin{remark}
\label{proposition_equivalence_A}
For the conditioned Schur complement ${\rho} \,|\, u$, the m.g.f.~from~(\ref{eq:non-central-chi-mgf}) of the mixture distribution in~(\ref{eq:non-central-chi-square-density}) and the p.g.f.~from~(\ref{equation_pgf_GA_s}) of the mixture weights in~(\ref{equation_new_alpha_k}) are related as
\begin{equation}
M_{{\rho} \,|\, u }(\theta) = (1-2\theta)^{-\nu/2} G_A\left(\frac{1}{1-2\theta}\right).
\label{eq:pgf-mgf-relation}
\end{equation}
\end{remark}

This reveals the equivalence between expressing the m.g.f.\ of a mixture of 
central $\chi^2$ distributions and expressing the
p.g.f.\ of the discrete random variable characterized by the mixture weights. 
Hereafter, we refer to the p.g.f.\ of the discrete random variable
as the p.g.f.\ of the mixture weights.

%Denote the p.d.f.~of $u$ by $h(u)$

%\begin{proposition}
%\label{prop:1}
%The p.d.f.~and m.g.f.~of the unconditioned $w_{11\cdot 2}/\sigma_{11\cdot 2}$ can be written as
%\begin{eqnarray}
%\label{eq:prop1}
%p(w) & = & \mathbb{E}\left[g(w;\nu, {\delta} = \lambda u) \,|\, u \right] = \int_0^1 g(w;\nu,\lambda u) \, h(u) \rm{d} u \nonumber \\
%& = & \sum_{k=0}^\infty \frac{(\lambda/2)^k}{k!}g_{\nu+2k}(w) \int_0^1 u^k e^{-\lambda u/2} h(u) \rm{d} u, \\
%\label{eq:prop2}
%M(\theta)& = &\frac{1}{(1-2\theta)^{\nu/2}} \int_0^1 \exp\left(\frac{\lambda\theta u}{1-2\theta}\right) h(u)\rm{d} u.
%\end{eqnarray}
%\end{proposition}

\subsection{General characterization of the unconditional distribution of the Schur complement}
\label{section_uncond_SC_preliminary}

Recall now from~(\ref{equation_definition_delta1}) that $\delta = \lambda u$, with $u$ being the random quadratic form defined in~(\ref{equation_definition_u_m1_Q2}).
Thus, we need to average the above distribution of ${\rho} \,|\, u$ over the distribution of $u$.
Letting the p.d.f.~of $u$ be denoted by $h(u)$, which is nonzero only for $u \in [0,1]$, the p.d.f.~and m.g.f.~of the scaled scalar Schur complement are given, respectively, by
\begin{align}
p_{{\rho}}(w)
=& \int_0^1 g(w;\nu,\lambda u) \, h(u) {\rm d}u \nonumber \\
=& \sum_{k=0}^\infty g_{\nu+2k}(w) \frac{(\lambda/2)^k}{k!} \int_0^1 u^k e^{-\lambda u/2} h(u) {\rm d}u,
\label{eq:prop1}
\end{align}
and 
\begin{equation}
\label{eq:prop2}
M_{{\rho}}(\theta) = (1-2\theta)^{-\nu/2} \int_0^1 \exp\left(\frac{\lambda\theta u}{1-2\theta}\right) h(u){\rm d}u.
\end{equation}

Thus, the study of the distribution of
${\rho}$ reduces to the study of the density $h(u)$ of $u$, which, in general, is difficult to derive.
%, which, in general, is difficult to deduce.
%in \eqref{eq:def-u}. 
%In general, $h(u)$ is difficult to derive. 
However, regardless of $h(u)$, it is straightforward to show that an m.g.f.--p.g.f.~relationship analogous to~(\ref{eq:pgf-mgf-relation}) holds also after removing the conditioning on $u$.
One can show this directly by averaging over $u$ in~(\ref{eq:pgf-mgf-relation}).

Nonetheless, for more details, let us arrive at the same result by casting the p.d.f.~from \eqref{eq:prop1} as the mixture of central $\chi^2$ p.d.f.s
\begin{equation}
\label{eq:alpha-k}
p_{{\rho}}(w) = \sum_{k=0}^\infty \beta_k \, g_{\nu + 2k}(w),
\end{equation}
with mixture weights
\begin{equation}
\label{eq:pf-of-N}
\beta_k = \frac{(\lambda/2)^k}{k!} \int_0^1 u^k e^{-\lambda u/2} h(u){\rm d}u,\ \ k=0,1,\ldots.
\end{equation}
We may view each weight $\beta_k$ as the probability that a nonnegative discrete random variable $B$ takes value $k$ because
\begin{align*}
%\label{eq:pf-of-N_sum}
\sum_{k=0}^{\infty} \beta_k &= \sum_{k=0}^\infty \beta_k \, \int_{0}^{\infty} g_{\nu + 2k}(w) \, {\rm d}w \\
&= \int_{0}^{\infty} \sum_{k=0}^\infty \beta_k \, g_{\nu + 2k}(w) \, {\rm d}w = \int_{0}^{\infty} p_{{\rho}}(w) \, {\rm d}w = 1,
\end{align*}
where the interchange of summation and integration can be justified by Fubini's theorem.  
Let $G_B(s)$ be the p.g.f.~of this random variable.

On the other hand, the m.g.f.~corresponding to the p.d.f.~in~(\ref{eq:alpha-k}) can be written as
\begin{align}
\label{eq:pgf-mgf-relation_temp}
M_{{\rho}}(\theta)
&= \int_0^\infty e^{\theta w} \sum_{k=0}^\infty \beta_k g_{\nu+2k}(w) \, {\rm d}w \nonumber \\
&= \sum_{k=0}^\infty \beta_k \int_0^\infty e^{\theta w} g_{\nu+2k}(w) \, {\rm d}w \nonumber \\
%&= \sum_{k=0}^\infty \beta_k (1-2\theta)^{-(\nu + 2k)/2} \\
&= (1-2\theta)^{-\nu/2} \sum_{k=0}^\infty (1-2\theta)^{-k} \beta_k.
\end{align}

%Now, by definition, the p.g.f.~of $N$ is given by
%\[
% G_N(s) = \mathbb{E}\left[s^N\right] = \sum_{k=0}^\infty s^k \beta_k.
%\]

\begin{proposition}
For the Schur complement $\rho$, regardless of $h(u)$, the mixture-distribution m.g.f. from (\ref{eq:pgf-mgf-relation_temp}) and the p.g.f.~of the corresponding mixture weights in~(\ref{eq:pf-of-N}) are related as follows:
\begin{equation}
M_{{\rho}}(\theta) = (1-2\theta)^{-\nu/2} G_B\left(\frac{1}{1-2\theta}\right).
\label{eq:pgf-mgf-relation2}
\end{equation}
\end{proposition}

%Notice that, the m.g.f.\ of noncentral $\chi^2$ distribution, i.e., 
%in \eqref{eq:non-central-chi-mgf}
%\[
% M_{{\rho} \,|\, u }(\theta) = \frac{1}{(1-2\theta)^{\nu/2}} e^{\frac{{\delta}}{2}\left(\frac{1}{1-2\theta}-1\right)},
%\]
%and the p.g.f.~of $Pois({\delta}/2)$ distributed mixture weights in~(\ref{eq:non-central-chi-square-density}), i.e.,
%\[
% G(s) = e^{({\delta}/2)(s-1)}
%\]
%are also related as in~(\ref{eq:pgf-mgf-relation}).
%, i.e., as
%\begin{eqnarray}
%M_{{\rho} \,|\, u }(\theta) = \frac{1}{(1-2\theta)^{\nu/2}}G\left(\frac{1}{1-2\theta}\right).
%\end{eqnarray}

%\footnote{As expected, considering the mixture weights in~(\ref{eq:non-central-chi-square-density}).} 

%\begin{remark}
%The above results reveal the equivalence between expressing the m.g.f.\ of a mixture of 
%central $\chi^2$ distributions and expressing the
%p.g.f.\ of the discrete random variable characterized by the mixture weights. 
%Therefore, we refer to the p.g.f.\ of the discrete random variable
%as the p.g.f.\ of the mixture weights.
%\end{remark}

The usefulness of this equivalence is supported below, as we characterize the distribution of $\rho$ for the case with $\rank\uOmega=1$.

\section{Distribution of the Schur complement with $\tilde \um_1 \neq \mzero$ and $\rank \uM = \rank \uOmega = 1$}
\label{sec:rank1}

The trivial case with $\tilde \um_1 = \mzero$, i.e., with $\um_1 = \uM_2 \uSigma_{22}^{-1} \usigma_{21}$, implies that $u = 0$, which, based on~(\ref{equation_SC_cond_distrib_u}), yields, simply, $\rho = {\rho} \,|\, u \sim \chi^2_{\nu } ({\delta})$, $\nu=n-p+1$,  \cite{siriteanu_twc_15}.
Therefore, hereafter, we consider the case with $\tilde \um_1 \neq \mzero$ and, for tractability, $\rank \uOmega=1$. 
Then, results shown above help recast more conveniently results deduced for $\uM_2=\mzero$ and $\uM_2\neq\mzero$ in \cite{siriteanu_twc_13} and \cite{siriteanu_twc_15}, respectively.

\subsection{Case with $\um_1 \neq \mzero$ and $\uM_2=\mzero$, i.e., $\rank \uM = \rank \uOmega = 1$}

Let the QR decomposition of $\uX_2 \in \mathbb{R}^{n \times (p-1)}$ be
\begin{equation*}
%\label{eq:QR1}
\uX_2 = \uH \uT,
\end{equation*}
where $\uT \in \mathbb{R}^{(p-1) \times (p-1)}$ is upper triangular with positive diagonal elements, and $\uH \in \mathbb{R}^{n \times (p-1)}$ satisfies $\uH'\uH = \mbI_{p-1}$.
Thus, for $\uM_2=\mzero$, $\uH$ has the uniform distribution over the Stiefel manifold $\{\uH \,|\, \uH'\uH=\mbI_{p-1}\}$.
Then, we can write for $\uQ_2$ the following
\begin{equation*}
%\label{eq:QR2}
\uQ_2 = \mbI_n - \uX_2 (\uX_2'\uX_2)^{-1} \uX_2' = \mbI_n - \uH \uH',
\end{equation*}
and for $u$ the following
\[
%\label{equation_u_intermsof_H}
u = \frac{\tilde \um_1' \uQ_2 \tilde \um_1}{\Vert \tilde \um_1\Vert^2} = \frac{\tilde \um_1'} {\Vert \tilde \um_1\Vert}
(\mbI_n - \uH \uH')\frac{\tilde \um_1}{\Vert \tilde \um_1\Vert} \sim Beta(\nu /2, (p-1)/2).
\]
%has the beta distribution with parameters $(\nu , p-1)$, 
Thus, $u$ has the central beta p.d.f.
\begin{equation}
\label{equation_hu_Beta}
h(u) = \frac{u^{\nu/2 -1} (1-u)^{(p-1)/2-1}}{B(\nu /2, (p-1)/2)}.
\end{equation}

With the p.d.f.~$h(u)$ from~(\ref{equation_hu_Beta}), the m.g.f.~of interest, i.e., $M_{{\rho}}(\theta)$ with ${\rho}=\lambda u$ from~(\ref{eq:prop2}), can be written, based on \cite[Eq.~(13.4.1)]{NIST_book_10}, as
\begin{align}
M_{{\rho}}(\theta)
&= (1-2\theta)^{-\nu/2} \int_0^1 \exp\left(\frac{\lambda\theta u}{1-2\theta}\right) h(u){\rm d}u \nonumber \\
&= (1-2\theta)^{-\nu/2} \hFoo \left(\frac{\nu }{2};\frac{n}{2}; \frac{\lambda\theta}{1-2\theta} \right).
\label{eq:1f1-mgf}
\end{align}
%\begin{eqnarray}
%\label{eq:1f1-mgf}
%M_{{\rho}}(\theta) = \frac{1}{(1-2\theta)^{\nu/2}} \hFoo \left(\frac{\nu }{2};\frac{n}{2}; \frac{\lambda\theta}{1-2\theta} \right),
%\end{eqnarray}
On the other hand, the positive weight for the corresponding mixture distribution characterized by \eqref{eq:alpha-k}, i.e., $\beta_k$ defined in~(\ref{eq:pf-of-N}), is given by \cite[Eq.~(13.4.1)]{NIST_book_10}
%\begin{align}
%\beta_k &= \frac{(\lambda/2)^k}{k!} \int_0^1 u^k e^{-\lambda u/2} h(u){\rm d} u 
%\nonumber \\
%& = 
%\frac{(\lambda/2)^k}{k!}
%\frac{B(k+\nu/2, (p-1)/2)}{B(\nu /2, (p-1)/2)}
% \hFoo \left(k+\frac{\nu }{2};k+\frac{n}{2}; -\frac{\lambda}{2} \right)
%\nonumber \\
%& = \frac{(\lambda/2)^k}{k!}
%\frac{{\Gamma}(k+\nu /2) {\Gamma}(n/2)}{{\Gamma}(\nu /2){\Gamma}(k+n/2)}
% \hFoo \left(k+\frac{\nu }{2};k+\frac{n}{2}; -\frac{\lambda}{2} \right).
%\label{eq:1f1-ak}
%\end{align}
\begin{align}
\beta_k
=& \int_0^1 \frac{u^{\nu /2 -1} (1-u)^{(p-1)/2-1} }{B(\nu /2+k, (p-1)/2)} \, \frac{(\lambda u/2)^k }{k!} e^{-\lambda u/2} {\rm d}u \nonumber \\
=& \frac{(\lambda/2)^k}{k!}
%\frac{B(k+\nu /2, (p-1)/2)}{B(\nu /2, (p-1)/2)} 
 \frac{(\nu /2)_k}{(n/2)_k}
 \hFoo \left(\frac{\nu }{2}+k;\frac{n}{2}+k; -\frac{\lambda}{2} \right).
\label{eq:1f1-ak}
\end{align}

%On the other hand, for the discrete random variable $B$ characterized by $\Pr(N = k)$ given by \eqref{eq:pf-of-N}, 
%we can write
%\begin{align*}
%%\label{eq:pf-of-N_detail}
%\Pr(N = k) = \beta_k
%=& \int_0^1 \frac{(\lambda/2)^k}{k!} e^{-\lambda u/2} u^k h(u){\rm d}u, \nonumber \\
%=& \int_0^1 \frac{u^{\nu /2 -1} (1-u)^{(p-1)/2-1} }{B((k + \nu )/2, (p-1)/2)} \, \frac{(\lambda u/2)^k }{k!} e^{-\lambda u/2} {\rm d}u,
%\end{align*}
%which coincides with (\ref{eq:1f1-ak}).

\medskip

\begin{remark}
Note from~(\ref{eq:1f1-ak}) that the distribution of random variable $B$ with $\Pr(B = k) = \beta_k$ is a beta mixture of Poisson distributions \cite[p.~42]{grandell-1997}, \cite{gurland-1958}, \cite{teunis-havelaar}, which, in turn, is a particular case of the
generalized hypergeometric factorial moment distribution \cite[Secs. 2.4.2, 6.1.13]{johnson-kemp-kotz}, \cite{kemp-kemp-1974}, \cite{lu-richards}. % \cite{sibuya-2006}, \cite{tripathi-gurland}.
\end{remark}

Because the integral and infinite series expressions deduced above for $\beta_k$ are involved, it is not trivial to deduce for $\beta_k$ the corresponding p.g.f.~$G_B(s) = \mathbb{E}\left[s^B\right] = \sum_{k=0}^\infty s^k \beta_k$.
On the other hand, by applying the general relationship revealed in \eqref{eq:pgf-mgf-relation2} for the m.g.f.~in~(\ref{eq:1f1-mgf}), the p.g.f.~for mixture weights $\beta_k$ can be expressed --- surprisingly simply --- as
\[
%\label{equation_GBs}
 G_B(s) = \hFoo \left(\frac{\nu }{2};\frac{n}{2}; \frac{\lambda}{2}(s-1) \right).
\]

\subsection{Case with $\um_1 \neq \mzero$ and $\uM_2 \neq \mzero$, but $\rank \uM = \rank \uOmega = 1$}

Recall from above that for $\um_1 \neq \mzero$ and $\uM_2=\mzero$, $u$ was shown to be central beta distributed.
On the other hand, when $\um_1 \neq \mzero$ and $\uM_2 \neq \mzero$, and, further, $\rank \uM = \rank \uOmega = 1$, we reveal next that $u$ is noncentral beta distributed (cf.\ \cite{hodges-1956}, \cite[Section 30.7.2]{johnson-kotz-balakrishnan-vol2}), \cite{seber-1963}.

\begin{proposition}
\label{prop:rank-1-los}
When $\um_1 \neq \mzero$, $\uM_2 \neq \mzero$, and $\rank \uM = \rank \uOmega = 1$, random variable $u$ has the distribution $Beta(\nu /2, (p-1)/2, \tau)$ with $\tau=\tr\left[\uSigma_{22}^{-1} \uM_2' \uM_2\right]$, i.e., $u$ has the p.d.f.
\begin{equation}
\label{eq:nc-beta-2}
h(u) = e^{-\tau/2} \sum_{l=0}^\infty \frac{(\tau/2)^l}{l!}
\frac{u^{\nu /2 - 1} (1-u)^{(p-1)/2 + l-1}}{B(\nu /2,(p-1)/2 + l)}.
\end{equation}
Then, the m.g.f.~of the scaled scalar Schur complement, the corresponding mixture weights, and their p.g.f.~can be expressed as follows:
%\begin{eqnarray}
%\label{eq:mgs-los-1}
%M_{{\rho}}(\theta) & = & \sum_{l=0}^\infty \frac{(\tau/2)^l}{l!}e^{-\tau/2} 
%\hFoo \left(\frac{\nu }{2};\frac{n}{2}+l; \frac{\lambda\theta}{1-2\theta}\right), \\
%\beta_k & = & \frac{(\lambda/2)^k}{k!}\sum_{l=0}^\infty \frac{(\tau/2)^l}{l!} 
%\frac{B(\nu /2+k,(p-1)/2+l)}{B(\nu /2,(p-1)/2+l)}
%\hFoo \left(k+\frac{\nu }{2};\frac{n}{2}+k+l; -\frac{\lambda}{2}\right) \nonumber \\
%& = & \frac{(\lambda/2)^k}{k!}\sum_{l=0}^\infty \frac{(\tau/2)^l}{l!} 
%\frac{{\Gamma}(\nu /2+k) {\Gamma}(n/2+l)}{{\Gamma}(\nu /2){\Gamma}(n/2+ k+l)} 
%\hFoo \left(\frac{\nu }{2}+k;\frac{n}{2}+k+l; -\frac{\lambda}{2} \right), \nonumber \\
%\label{equation_pgf_Gs}
%G(s)& = & \sum_{l=0}^\infty \frac{(\tau/2)^l}{l!}e^{-\tau/2} 
%\hFoo \left(\frac{\nu }{2};\frac{n}{2}+l; \frac{\lambda}{2}(s-1) \right).
%\end{eqnarray}
\begin{align}
M_{{\rho}}(\theta)
&= (1-2\theta)^{-\nu/2} e^{-\tau/2} \sum_{l=0}^\infty \frac{(\tau/2)^l}{l!}
\hFoo \left(\frac{\nu }{2};\frac{n}{2}+l; \frac{\lambda\theta}{1-2\theta}\right) \nonumber \\
&= (1-2\theta)^{-\nu/2} e^{-\tau/2}
\Phi_2\left(\frac{\nu }{2},\frac{n}{2};\frac{n}{2};\frac{\lambda\theta}{1-2\theta},\frac{\tau}{2}\right),
\label{eq:mgs-los-1}
\end{align}
\begin{align}
\beta_k
&= \frac{(\lambda/2)^k}{k!}\sum_{l=0}^\infty \frac{(\tau/2)^l}{l!} 
%\frac{B(\nu /2+k,(p-1)/2+l)}{B(\nu /2,(p-1)/2+l)}
% \nonumber \\ & \quad \times 
\frac{(\nu /2)_k}{(n/2+l)_k}
\hFoo \left(\frac{\nu }{2}+k;\frac{n}{2}+k+l; -\frac{\lambda}{2}\right) \nonumber \\
&= \frac{(\lambda/2)^k}{k!} \frac{(\nu /2)_k}{(n/2)_k}
\Phi_2\left(\frac{\nu }{2}+k,\frac{n}{2};\frac{n}{2}+k;-\frac{\lambda}{2},\frac{\tau}{2}\right),
\label{equation_alpha_k_prop}
\end{align}
\begin{align}
G(s)
&= e^{-\tau/2} \sum_{l=0}^\infty \frac{(\tau/2)^l}{l!}
\hFoo \left(\frac{\nu }{2};\frac{n}{2}+l; \frac{\lambda}{2}(s-1) \right) \nonumber \\
&= e^{-\tau/2} \Phi_2\left(\frac{\nu }{2},\frac{n}{2};\frac{n}{2};\frac{\lambda}{2}(s-1),\frac{\tau}{2}\right).
\label{equation_pgf_Gs}
\end{align}
\end{proposition}

\medskip

We remark that, to the best of our knowledge, the p.g.f.~in~(\ref{equation_pgf_Gs}) has not been studied before now.

%Let us present the consequences of this proposition. 
%From \eqref{eq:nc-beta-2} we have the following expressions for the m.g.f., weights and the p.g.f.
%\begin{equation}
%\label{eq:mgs-los-1}
%M(\theta) = \sum_{l=0}^\infty \frac{(\tau/2)^l}{l!}e^{-\tau/2} 
%\hFoo \left(\frac{\nu }{2};\frac{n}{2}+l; \frac{\lambda\theta}{1-2\theta}\right).
%\end{equation}
%Thus, the weight $\beta_k$ is given by
%\begin{align*}
%\beta_k &= \frac{(\lambda/2)^k}{k!}\sum_{l=0}^\infty \frac{(\tau/2)^l}{l!} 
%\frac{B(\nu /2+k,(p-1)/2+l)}{B(\nu /2,(p-1)/2+l)}
%\hFoo \left(k+\frac{\nu }{2};\frac{n}{2}+k+l; -\frac{\lambda}{2}\right) \\
%&=\frac{(\lambda/2)^k}{k!}\sum_{l=0}^\infty \frac{(\tau/2)^l}{l!} 
%\frac{{\Gamma}(\nu /2+k) {\Gamma}(n/2+l)}{{\Gamma}(\nu /2){\Gamma}(n/2+ k+l)} 
%\hFoo \left(\frac{\nu }{2}+k;\frac{n}{2}+k+l; -\frac{\lambda}{2} \right)
%%\qquad(\text{check!})
%\end{align*}
%with p.g.f.
%\begin{eqnarray}
%G(s)=\sum_{l=0}^\infty \frac{(\tau/2)^l}{l!}e^{-\tau/2} 
%\hFoo \left(\frac{\nu }{2};\frac{n}{2}+l; \frac{\lambda}{2}(s-1) \right),
%\end{eqnarray}
%which does not seem to have been studied in literature.

\begin{proof}
Let us provide for~(\ref{eq:nc-beta-2}) a new proof,
%\footnote{Different from the proof in \cite{siriteanu_twc_15}.}
different from the one in \cite{siriteanu_twc_15},
based on known distributional results on Hotelling's $T^2$-statistic (cf.\ \cite[Ch.~5]{anderson-3rd}).

First, let us transform the parameter matrices $\uSigma$ and $\uM$ to %canonical form.
simpler forms.
This can be done by noting that $\tilde \ux_1$ had been obtained in \eqref{eq:x1-sub-X2} by the following column-wise linear transformation of $\uX$:
\[
\tilde \ux_1 = \ux_1 - \uX_2 \uSigma_{22}^{-1} \usigma_{21}
 = \begin{pmatrix} \ux_1 & \uX_2 \end{pmatrix}
 \begin{pmatrix} 1 \\ -\uSigma_{22}^{-1} \usigma_{21} \end{pmatrix}
 = \uX \begin{pmatrix} 1 \\ -\uSigma_{22}^{-1} \usigma_{21}\end{pmatrix}.
\]
%Applying the same transformation to $\uSigma$ yields
%\[
%\uSigma=\begin{pmatrix} \sigma_{11} & \usigma_{21}' \\ \usigma_{21} & \uSigma_{22} \end{pmatrix}
% \begin{pmatrix} 1 \\ - \uSigma_{22}^{-1} \usigma_{21} \end{pmatrix}
% = \begin{pmatrix} \sigma_{11} - \usigma_{21}' \uSigma_{22}^{-1} \usigma_{21} \\ \mzero \end{pmatrix}
% = \begin{pmatrix} \sigma_{11\cdot 2} \\ \mzero
%\end{pmatrix},
%\]
%%This also transforms the first row of the covariance matrix $\uSigma$ from~(\ref{equation_covariance}), i.e., $(\sigma_{11}, \usigma_{21}')$, into $(\sigma_{11\cdot 2},0)$, 
%with Schur complement $\sigma_{11\cdot 2}$ appearing in~(\ref{equation_definition_rho1}) as a scale factor for Schur complement $w_{11\cdot 2}$.
%% considering $\tilde \ux_1/\sqrt{\sigma_{11\cdot 2}}$ instead of $\tilde \ux_1$ 
%% we can assume $\sigma_{11\cdot 2}=1$. 
%% Now write
%% \begin{equation}
%% \label{eq:sigma22-half}
%% \uSigma_{22} = (\uSigma_{22}^{1/2})'\uSigma_{22}^{1/2}
%% \end{equation}
%% and transform $\uX_2$ to $\tilde \uX_2 =\uX_2 \uSigma_{22}^{-1/2}$. 
%% Then $\uSigma_{22}$ is transformed to $\mbI_{p-1}$, but
%% $\uX_2$ and $\tilde \uX_2$ span the same columns space and the orthogonal projector is unchanged
%% (i.e. $\uX_2 (\uX_2'\uX_2)^{-1} \uX_2' = \tilde \uX_2 (\tilde \uX_2'\tilde \uX_2)^{-1} \tilde \uX_2'$).
%% Hence the distribution of $w_{11\cdot 2}$ does not change by this transformation.
%% Summarizing the above column transformations, by
Thus, employing
\[
%\uA = \begin{pmatrix} \frac{1}{\sqrt{\sigma_{11\cdot 2}}} & 0 \\ 
%  - \frac{\uSigma_{22}^{-1} \usigma_{21} }{\sqrt{\sigma_{11\cdot 2}}} & \mbI_{p-1} \end{pmatrix},
\uA = \begin{pmatrix} 1 & \mzero' \\ -\uSigma_{22}^{-1} \usigma_{21} & \mbI_{p-1} \end{pmatrix}
\]
for the column-wise linear transformation $\uX \mapsto \uX \uA$ yields
%% \begin{eqnarray}
%% X \mapsto (\tilde \ux_1/\sqrt{\sigma_{11\cdot 2}}, \uX_2 \uSigma_{22}^{-1/2})
%% =X A, \qquad A=\begin{pmatrix} 1/\sqrt{\sigma_{11\cdot 2}} & 0 \\ 
%%     - \uSigma_{22}^{-1} \usigma_{21}/\sqrt{\sigma_{11\cdot 2}} & \uSigma_{22}^{-1/2}
%%   \end{pmatrix},
%% \end{eqnarray}
%\begin{eqnarray}
% \uX & \mapsto & \uX \uA 
% = \begin{pmatrix} \ux_1 & \uX_2
%\end{pmatrix} \begin{pmatrix} \frac{1}{\sqrt{\sigma_{11\cdot 2}}} & 0 \\
%  - \frac{\uSigma_{22}^{-1} \usigma_{21} }{\sqrt{\sigma_{11\cdot 2}}} & \mbI_{p-1}
%   \end{pmatrix} = \begin{pmatrix}\frac{\tilde \ux_1}{\sqrt{\sigma_{11\cdot 2}}} & \uX_2 \end{pmatrix}, \nonumber \\
% \label{eq:m-ma}
%\uM & \mapsto & \uM \uA = \begin{pmatrix} \um_1 & \uM_2
%\end{pmatrix} \begin{pmatrix} \frac{1}{\sqrt{\sigma_{11\cdot 2}}} & 0 \\ 
%   - \frac{\uSigma_{22}^{-1} \usigma_{21} }{\sqrt{\sigma_{11\cdot 2}}} & \mbI_{p-1}
%  \end{pmatrix} = %(\tilde \um_1/\sqrt{\uSigma_{11.2}}, \uM_2 \uSigma_{22}^{-1/2}).
%\begin{pmatrix}\frac{\tilde \um_1}{\sqrt{\sigma_{11\cdot 2}}}& \uM_2\end{pmatrix},
%\end{eqnarray}
\begin{equation}
 \uX \, \mapsto \, \uX \uA 
 = \begin{pmatrix} \ux_1 & \uX_2 \end{pmatrix}
 \begin{pmatrix} 1 & \mzero' \\ -\uSigma_{22}^{-1} \usigma_{21} & \mbI_{p-1} \end{pmatrix}
 = \begin{pmatrix} \tilde \ux_1 & \uX_2 \end{pmatrix},
\label{eq:m-ma}
\end{equation}
which preserves the Schur complement $w_{11\cdot 2}$ of $\uW=\uX'\uX$.
On the other hand, this transformation induces the transformation on the parameter space of $\uM$ and $\uSigma$ as follows:
\begin{align*}
\uM \, &\mapsto \, \uM \uA
% = \begin{pmatrix} \um_1 & \uM_2 \end{pmatrix}
% \begin{pmatrix} 1 & 0 \\ -\uSigma_{22}^{-1} \usigma_{21} & \mbI_{p-1} \end{pmatrix}
 = \begin{pmatrix} \tilde \um_1 & \uM_2\end{pmatrix}, \\
\uSigma \, &\mapsto \, \uA' \uSigma \uA
 = \begin{pmatrix} \sigma_{11\cdot 2} & \mzero' \\ \mzero & \uSigma_{22} \end{pmatrix}.
\end{align*}
%preserving $\uX_2$ and $\uM_2$, and yielding $\rank \uM = \rank (\uM \uA) = 1$.
Note that this preserves $\rank \uM = \rank (\uM \uA) = 1$.
%$\uSigma$ is transformed to $I_p$. 
%Then, $\uM= \mathbb{E}[X]$ undergoes the same transformation, i.e.,
%\begin{equation}
%\label{eq:m-ma}
%\uM \mapsto \uM \uA = %(\tilde \um_1/\sqrt{\uSigma_{11.2}}, \uM_2 \uSigma_{22}^{-1/2}).
%\left(\frac{\tilde \um_1}{\sqrt{\sigma_{11\cdot 2}}}, \uM_2\right),
%\end{equation}
% Also note that in \eqref{eq:sigma22-half}
% we can replace $\uSigma_{22}^{1/2}$ by $Q\uSigma_{22}^{1/2}$ where $Q$
% is orthogonal.  This further orthogonal transformation will be used to transform $\uM$.

Second, let us apply the row-wise transformation $\uX \mapsto \uP \uX$, with orthogonal matrix $\uP \in \mathbb{R}^{n \times n}$ specified below. 
Then, $\uW=\uX'\uX$ and the covariance structure $\uSigma$ of $\uX$ are unchanged, whereas mean matrix is transformed as $\uM \mapsto \uP \uM$.

Since $\rank (\uM \uA)=1$, there exists an $n$-dimensional column vector $\ua$ and a $p$-dimensional column vector $\ub=(b_1,\dots, b_p)'$ such that
\[
 \uM \uA=\ua\ub',
\]
where we may assume $\Vert \ua \Vert=1$, because %either $\ua$ or $\ub$ can be assumed not to be zero 
$\ua$ and $\ub$ are nonzero vectors
and $\ua\ub'=(\ua/c) c\ub'$ for any constant $c \neq 0$.
Then, by \eqref{eq:m-ma} and~(\ref{equation_definition_lambda}) we have 
%\begin{eqnarray}
%\ua=\pm \tilde \um_1 / \Vert \tilde \um_1\Vert, \quad \text{and} \quad \lambda = \frac{\Vert \tilde \um_1\Vert^2}{\sigma_{11\cdot 2}} = b_1^2, b_1 = \frac{\Vert \tilde \um_1\Vert}{\sqrt{\sigma_{11\cdot 2}}} = \sqrt{\lambda}.
%\end{eqnarray}
\[
 \ua=\frac{\tilde \um_1}{\Vert \tilde \um_1\Vert}, \quad \text{and} \quad
% b_1 = \frac{\Vert \tilde \um_1\Vert}{\sqrt{\sigma_{11\cdot 2}}} = \sqrt{\lambda}.
 b_1 = \Vert\tilde \um_1\Vert = \sqrt{\lambda \sigma_{11\cdot 2}}.
\]
Constructing $\uP$ so that its first row is $\ua'$ and its remaining rows are orthogonal to $\ua$ yields
\begin{equation}
\uP \, \uM \, \uA = \uP \, \ua \, \ub'
 = \begin{pmatrix} 1 \\ 0 \\ \vdots \\ 0 \end{pmatrix} \ub'
 = \begin{pmatrix} \sqrt{\lambda \sigma_{11\cdot 2}} & b_2 & b_3 & \dots & b_p \\
 0      & 0      & 0      & \dots & 0 \\
 \vdots & \vdots & \vdots &       & \vdots \\
 0      &      0 &   0    & \dots & 0
\end{pmatrix}.
\label{equation_PMAQ0}
\end{equation}
The assumption $\uM_2\neq \mzero$ implies that $\ub_2=(b_2,\dots,b_p)'\neq \mzero$.
Therefore, it suffices to prove the proposition when $\uSigma$ and $\uM$ are as follows:
\[
% \uSigma=\begin{pmatrix} 1 & 0 \\ 0 & \uSigma_{22} \end{pmatrix}, \ 
 \uSigma = \begin{pmatrix} \sigma_{11\cdot 2} & \mzero' \\ \mzero & \uSigma_{22} \end{pmatrix}, \quad 
% I_p, 
% \ \uM=\begin{pmatrix} 1 & 0 \\ 0 & \uSigma_{22} 
% \end{pmatrix}
% I_p, 
% \ \ua \, \ub = 
% \begin{pmatrix} \sqrt{\lambda} & \sqrt{\tau} & 0 & \dots & 0 \\
%  0 & 0 & 0 & \dots & 0 \\
% \vdots & \vdots & \vdots& \vdots & \vdots \\
%  0    &      0 &   0 & \dots & 0 \\
% \end{pmatrix}. 
\uM =
%\begin{pmatrix} \sqrt{\lambda} & b_2 & \dots & b_p \\
\begin{pmatrix} \sqrt{\lambda \sigma_{11\cdot 2}} & b_2 & \cdots & b_p \\
 0     & 0      & \cdots & 0 \\
\vdots & \vdots &        & \vdots \\
 0     & 0      & \cdots & 0
\end{pmatrix}. 
\]

Note that in the matrix in~(\ref{equation_PMAQ0}), in the first column, which corresponds to $\tilde \um_1$ from~(\ref{eq:x1givenx2}), only the first element is nonzero.
Therefore, the random variable $u$ from~(\ref{equation_definition_u_m1_Q2}) is the $(1,1)$ element of the matrix
%\begin{eqnarray}
$ \uQ_2 = \mbI_n -\uX_2 (\uX_2'\uX_2)^{-1} \uX_2'$.
%= \mbI_n - HH'.
%\end{eqnarray}
% Denoting the first row of $\uH$ as $(h_{11},\dots, h_{1,p-1})$ we have
%Denoting the first row of $\uX_2$ $(x_{12},\dots, x_{1p})$ we have
This element can be written, by partitioning $\uX_2$ into its first row and its remaining rows as
\[
 \uX_2 = \begin{pmatrix} \ux_{12}' \\ \uX_{22} \end{pmatrix},
\]
and by applying the matrix inversion lemma, as follows:
\begin{align*}
u
&= 1 - \ux_{12}' (\uX_{22}' \uX_{22} + \ux_{12} \ux_{12}')^{-1} \ux_{12} \\
&= 1 - \ux_{12}' \left((\uX_{22}' \uX_{22})^{-1} - \frac{(\uX_{22}' \uX_{22})^{-1} \ux_{12} \ux_{12}'(\uX_{22}' \uX_{22})^{-1}}{1+\ux_{12}' (\uX_{22}' \uX_{22})^{-1}\ux_{12}}\right) \ux_{12} \\
&= 1 - \frac{\ux_{12} (\uX_{22}' \uX_{22})^{-1}\ux_{12}'}{1+\ux_{12}' (\uX_{22}' \uX_{22})^{-1}\ux_{12}}
 = \frac{1}{1+\ux_{12}' (\uX_{22}' \uX_{22})^{-1}\ux_{12}}.
\end{align*}
Here $\ux_{12}'$ and the rows of $\uX_{22}$ are independent Gaussian vectors with means $\ub_2$ and $\mzero$, respectively, but the same covariance matrix $\uSigma_{22}$. Then, by \cite[Theorems~5.2.2,~5.4.1]{anderson-3rd}, the scalar random variable
$ T^2 =\ux_{12}' (\uX_{22}' \uX_{22})^{-1}\ux_{12}$ is a statistic of Hotelling's $T^2$-type, and therefore has the noncentral $F$-distribution with p.d.f.
\[
f_{T^2 }(t) = e^{-\tau/2} \sum_{l=0}^\infty \frac{(\tau/2)^l}{l!}
 \frac{1}{B(\nu /2,(p-1)/2 + l)} \frac{t^{p/2-1+l}}{(1+t)^{n/2+l}},
\]
where $\tau$ is the trace of the noncentrality matrix for the Wishart distribution of $\uW_{22}$, i.e.,
\[
%\tr \uM_2 '\uSigma_{22}^{-1/2} (\uSigma_{22}^{-1/2})' \uM_2' = \tr \uM_2 \uSigma_{22}^{-1} \uM_2'
\tau= \tr\left[\uSigma_{22}^{-1} \uM_2' \uM_2\right] = \tr\left[\uP \uM_2 \uSigma_{22}^{-1} (\uP \uM_2)'\right] = \ub_2\uSigma_{22}^{-1}\ub_2'.
\]
%Because $\ux_{12}' (\uX_{22}' \uX_{22})^{-1}\ux_{12} = 1 - \frac{1}{u} $, a simple transformation 
Then the transformation $u=1/(1+T^2 )$ yields for $u$ the p.d.f.~from \eqref{eq:nc-beta-2}.

Finally, using this p.d.f.~along with results from Section~\ref{sec:pdf} produces~(\ref{eq:mgs-los-1}),~(\ref{equation_alpha_k_prop}), and~(\ref{equation_pgf_Gs}).
\end{proof}

% \begin{remark}
% In \cite{siriteanu_twc_15},
% the series expansion apparently involved negative weights.
% But when we add these weights for each degree of freedom of central chi-square density, then the weights become positive and add up to 1.
% Note that weight is uniquely determined.
% \end{remark}

%\footnotesize
%\small

%\bigskip
%\bigskip

\section*{Acknowledgments}
A. Takemura acknowledges the support of the JSPS Grant-in-Aid for Scientific Research No. 25220001.

\begin{center}
\large{References}
\end{center}

\bibliographystyle{abbrv}
%\bibliographystyle{elsarticle-harv} 
%\bibliographystyle{elsarticle-num} 
%\bibliography{ghfd-wishart}

\end{document}